\documentclass[preprint,12pt]{elsarticle}

\usepackage{amscd,amsfonts,amssymb,amsmath,amsthm,latexsym,array,hhline,xcolor,graphicx}

\usepackage{latexsym}
\usepackage{times}

\newcommand\F{\mbox{I\kern-2pt F}}

\newcommand\cE{{\cal E}}

\newcommand\cF{{\cal F}}
\newcommand\cG{{\cal G}}

\newcommand\e{{\varepsilon}}

\newcommand{\RNumb}[1]{\uppercase\expandafter{\romannumeral #1\relax}}

\def\E{{\bf E}}
\def\P{{\bf P}}
\def\Chi{{\bf 1}}

\def\bbr{{\mathbb R}}

\newcommand\fdem{$\Box$}
\newcommand\beq{\begin{equation}}
\newcommand\eeq{\end{equation}}
\newcommand\bea{\begin{eqnarray}}
\newcommand\eea{\end{eqnarray}}
\newcommand\bean{\begin{eqnarray*}}
\newcommand\eean{\end{eqnarray*}}

\newtheorem{theorem}{Theorem}[section]
\newtheorem{corollary}[theorem]{Corollary}      %[section]
\newtheorem{lemma}[theorem]{Lemma}              %[section]
\theoremstyle{definition}

\usepackage[pdftex,colorlinks,urlcolor=red,citecolor=blue,linkcolor=red]{hyperref}

\begin{document}

\begin{frontmatter}

\title{Ruin Probabilities for a Sparre Andersen Model with Investments}

\author{Ernst Eberlein}
\address{University of Freiburg, Department of Mathematical Stochastics, Ernst-Zermelo-Str. 1, 79104 Freiburg, Germany}
\ead{eberlein@stochastik.uni-freiburg.de}

\author{Yuri Kabanov}
\address{Lomonosov Moscow State University, Russia,  and 
Laboratoire de Math\'ematiques, Universit\'e de Franche-Comt\'e, 16 Route de Gray,
25030 Besan\c{c}on, cedex, France}
\ead{ykabanov@univ-fcomte.fr}

\author{Thorsten Schmidt}
\address{University of Freiburg, Department of Mathematical Stochastics, Ernst-Zermelo-Str. 1, 79104 Freiburg, Germany}
\ead{thorsten.schmidt@stochastik.uni-freiburg.de}
\tnotetext[label1]{This work was done during the stay of the  second author  at the Freiburg Institute for Advanced Studies (FRIAS), Germany; it was supported 
by the Russian Science Foundation grant  
20-68-47030.
The third author gratefully acknowledges the support by the DFG with grant SCHM 2160/15-1.}

\date{\today}

\begin{abstract}
We study a Sparre Andersen model in which the business activity of the company is described by a compound renewal process with drift 
assuming that the capital reserves are invested in a risky asset. The price of the latter is assumed to evolve according to a geometric L\'evy process. 
We prove that the asymptotic behavior of the ruin probability  depends to a large extent only on the properties of the price process.
\end{abstract}


%  \subclass{60G44}
%  \medskip
% \noindent
%  {\bf JEL Classification} G22 $\cdot$ G23

\end{frontmatter}

%\newpage
%\acknowledgement{This work was done during the stay of the  second author  at FRIAS.} 

%\titlerunning{Sparre Andersen model with Investments}

\section{Introduction}  
In the classical  collective risk theory  the Sparre Andersen model (called also renewal risk model) is a generalisation of the  Lundberg--Cram\'er model in which it is assumed that the claims count is a renewal process rather than a Poisson process, i.e. the 
times between consecutive claims are independent random variables which are not necessarily exponentially distributed, see, e.g., Grandell's  book \cite{Gr}.   
The word ``classical" denotes the situation in which the reserves of the insurance company are kept in a bank account which does not earn any interest. 
In the modern literature  more realistic models are considered, namely those in which the capital is invested, in full or partially, in a risky financial asset. It was shown in \cite{FrKP} that the financial risk leads to a dramatic change in the behavior of the ruin probability $\Psi(u)$ as a function of the initial capital $u$: instead of its comfortable exponential decrease with the growth of the available capital, in case of the investment of the reserves in the stock market, the ruin probability $\Psi(u)$  decreases as a power function. Moreover, the ruin is imminent, $\Psi\equiv 1$, if the volatility of the stock market is too high. 
Because of the practical importance of this phenomenon, actuarial models with investments attracted a lot of attention. 
Most of the papers deal with asymptotic studies in settings stemming  from the  Lundberg--Cram\'er model and can be included in the beautiful mathematical framework of the generalised Ornstein--Uhlenbeck process suggested and studied by Paulsen, 
see \cite{Paul-93} -- \cite{Paul-G}. Paulsen used the powerful method of distributional equations. A detailed analysis of the available results is beyond the scope of the present note and we refer the reader to the recent paper \cite{KP2020}. 

 In contrast to the Lundberg--Cram\'er line of research there are very few results on ruin probabilities for Sparre Andersen type models with investments. The paper by Albrecher, Constantinescu, and Thomann \cite{ACT} is one of the rare studies to the latter, see also references therein. Their approach is heavily based on the techniques of continuous time Markov processes,   
 
In the present note we consider a Sparre Andersen type non-life insurance
model with investment in a risky asset. The price of this asset evolves according to an exponential L\'evy process, while the business activity of the company is described by a compound renewal process with positive drift and negative jumps. 
The asymptotic result which is achieved is basically of the same form as the main result of \cite{KP2020}. It requires only very weak assumptions on the interarrival times, like existence of the exponential moments, which are fulfilled in the case of the exponential distribution.  As in \cite{KP2020}, the main assumption is that
the cumulant generating function $H:\gamma\mapsto \ln \E e^{-\gamma V_1}$ of the log price process $V$ has a strictly positive root $\beta$
which does not lie on the boundary of the effective domain of $H$.
Assuming that the claims are random variables, whose power of order $\beta$ is integrable, the ruin probability decays with the rate $u^{-\beta}$  (see Theorem \ref{main}). 
 %Recall that in the case where the price process is a geometric
 %Brownian motion $\beta = 2a/\sigma^2-1$.       

As in \cite{KP2020}, our proofs  are based on  recent progress in the theory of distributional equations.

\section{The model}
\label{model}
Let  $(\Omega,\cF,{\bf F}=(\cF_t)_{t\ge 0},\P)$ be a stochastic basis on which there is given a non-deterministic 
L\'evy process $R=(R_t)_{t\ge 0}$ with  $\Delta R >-1$ and an independent 
compound renewal process $P=(P_t)_{t \ge 0}$ with drift $c>0$ and negative jumps. We denote by $T_n$, $n \ge 1$, the successive jump instants of the latter process and put $F(t):=\P(T_1\le t)$. 
We associate with $R$ the stochastic exponential 
$\cE(R)$ which is, due to the assumption on the jumps, a strictly positive process. 
%We denote $(a,\sigma^2,\Pi)$  the L\'evy triplet of $R$.  

% We denote $p_P(dt,dx)$ the jump measure of the latter with  its mean measure $\Pi_P(dx)dt$ where $\Pi_P(\R)<\infty$.  We assume that  $\Pi_P(\R_+)>0$ and $\Pi_P(\R_-)>0$ (some comments for  the cases where  $\Pi_P$ charges only of the half-axes will be also given). 

We will study the process $X=X^{u}$, $u>0$,  which is defined as the solution of the non-homogeneous linear stochastic 
%differential 
equation 
\beq
 \label{risk}
 X_t=u+  \int_0^t  X_{s-}  dR_s +P_t. 
 \eeq
%where $R_t$ is the relative price  process of the risky asset and $P$ is a compound Poisson process with drift $c\in \R$ representing the business activity of the insurance company, 
%$u>0$ is the initial capital at time zero; we assume that  $\sigma^2>0$. 
The process $X$ can be written in the ``dot" notation of stochastic calculus for  as $X=u+X_{-}\cdot R+P$ or, in  
its ``differential" form, as $dX_t= X_{t-}  dR_t +dP_t$, $X_0=u$. %We denote $T_n$ the successive instants of  jumps  of the Poisson process $N_t=p_P([0,t],\R)$.  
The compound renewal  process $P$   is usually represented in the form 
\beq
\label{Pt}
P_t=ct+\sum_{i=1}^{N_t}\xi_i,  
\eeq 
where $N$ is a counting renewal process with interarrival times $T_i-T_{i-1}$, $i\ge 1$, which form an i.i.d. sequence, independent of the i.i.d. sequence $\xi_i=\Delta P_{T_i}$, $i \ge 1$.
The case considered here ($c>0$ and $\xi_i<0$) corresponds to the non-life insurance  situation. For $R= 0$ we have the classical Sparre Andersen model (which we already excluded by the assumptions made).  

In the actuarial context  $X=X^u$ represents  the dynamics of the reserves of an insurance company that invests the total of these reserves  in a stock with the price process 
$S=\cE(R)$ which has unit as its initial value.
Having in mind that $dR=dS/S_-$, it is natural to call $R$ {\it relative price process}, while $V:=\ln \cE(R)$ is usually referred to as the {\it log price process}.
Thus the price process $S$ which we introduced here as a stochastic exponential can alternatively be expressed as an ordinary exponential $S=\exp(V)$. For the interplay between stochastic and ordinary exponentials see Theorem 3.49 in \cite{EK}. Note that $V$ is again a L\'evy process. In general this process will have jumps of arbitrary size in both directions. Typical examples which are used in financial models are generalized hyperbolic processes and various subclasses such as hyperbolic, normal inverse Gaussian or variance gamma processes.

 %B(t,T)= \exp\bigg(-\int_t^T f(t,s)ds\bigg).

The processes $R$ and $P$ generate  the filtration ${\bf F}^{R,P}=(\cF^{R,P}_t)_{t\ge 0}$.  
 Define  $\tau^u:=\inf \{t: X^u_t\le 0\}$ (the instant of ruin) and
 $\Psi(u):=\P(\tau^u<\infty)$ (the ruin probability). %Consequently 
 %$\Phi(u):=1-\Psi (u)$
 %(the survival probability).

\smallskip
%Note that the same formula  (\ref{Pt}) $N$ can be an independent renewal process, that is, the counting process in which   
%the lengths of the interarrival intervals $T_n-T_{n-1}$ form an i.i.d. sequence. 
%In the collective risk theory this corresponds to the Sparre Andersen model. 
 
%\section{Preliminaries from  the theory of L\'evy processes}
%\label{prel}

Let  $(a,\sigma^2,\Pi)$  be the L\'evy triplet 
of $R$  corresponding to the standard truncation function
$h(x):=xI_{\{|x|\le 1\}}$.  Putting $\bar h(x):=xI_{\{|x|> 1\}}$, we can write 
the canonical decomposition of $R$ in the form
\beq
\label{Rt}
R_t=at+\sigma W_t+h*(\mu -\nu)_t+\bar h*\mu_t. 
\eeq
Here $W$ is a standard Wiener process,  $\mu(dt,dx)$ is the random measure of jumps  of $R$,  i.e.  the Poisson random measure  with  compensator   
$\nu(dt,dx)=dt\Pi (dx)$. 
%The alternative L\'evy process $V$ (log price) would then be given by
%$$
%V_t=R_t-\frac 12 \sigma^2 t + (\ln (1+x)-x)*\mu_t.
%$$

For notations and basic facts we recommend the books \cite{JS}, Ch. 2,  and \cite{EK}, Chs. 2 and 4. 
As in \cite{KP2020} the symbols $\Pi(f)$ or  $\Pi(f(x))$ stand for the integral of $f$ with respect to the L\'evy measure $\Pi$. Recall that $\Pi(x^2\wedge 1)<\infty$.

Due to the independence of $R$ and $P$, these processes have no common jumps. The solution of the generalised stochastic exponential (\ref{risk})
can be represented (see \cite{EK}, Proposition 3.48) as
\beq
\label{uY}
X_t^u:=\cE_t(R)(u+\cE^{-1}_-(R)\cdot P_t). 
\eeq 

%\smallskip
%\noindent
%{\bf Notations:} 
%As in \cite{JS}, we use $*$ for the standard notation of stochastic calculus  for integrals with respect to random measures. For instance,    
%$$
%h*(\mu -\nu)_t=\int_0^t\int h(x)(\mu-\nu)(ds,dx). 
%$$
%We hope that the reader will be not confused that $f(x)$  may denote the whole function $f$ or its value at $x$; the typical example is $\ln (1+x)$ explaining why such a flexibility is convenient. 
%The symbols $\Pi(f)$ or  $\Pi(f(x))$ stands for the integral of $f$ with respect to the measure $\Pi$. 
%Recall that
% $$
%%\Pi(h^2)+\Pi(I_{\{|x|> 1\}})=  
%\Pi(x^2\wedge 1):=\int (x^2\wedge 1) \Pi(dx)<\infty
% $$
%and the condition $\sigma=0$ and $\Pi(|h|)<\infty$ is necessary and sufficient for  $R$ to have trajectories of (locally) finite variation, see Prop. 3.9 in \cite{CT}.  

Since we assume that $\Delta R >-1$ and $R$ is nondeterministic, the  L\'evy measure $\Pi$ is concentrated on the interval $]-1,\infty[$ and the L\'evy characteristics  $\sigma^2$ and $\Pi$ do not vanish simultaneously.

\medskip
Recall that the stochastic exponential $\cE_t(R)$ has explicit form 
$$
\cE_t(R)=e^{R_t -\frac 12\sigma^2t+\sum_{s\le t}(\ln (1+\Delta R_s)-\Delta R_s)}. 
$$
 The log price $V =\ln \cE(R)$ can then be expressed as
\beq
\label{V}
V_t=at-\frac 12 \sigma^2 t + \sigma W_t+ h*(\mu-\nu)_t+(\ln (1+x)-h)*\mu_t. 
\eeq
Its L\'evy triplet is  $(a_V,\sigma^2,\Pi_V)$ where 
$$
a_V=a-\frac{\sigma^2}{2}+\Pi(h(\ln(1+x))-h)  
$$ and 
$\Pi_V=\Pi\varphi^{-1}$,  $\varphi: x\mapsto \ln (1+x)$.  

\smallskip
The cumulant generating function $H:q\to \ln  \E\,e^{-qV_1}$ of the random variable $V_1$ admits an explicit expression 
\beq
\label{sec.Main.1}
H(q):=-a_V q+\frac{\sigma^{2}}{2}q^{2}
+\Pi \big(e^{-q \ln(1+x)}-1+q h(\ln (1+x))\big).  
\eeq
Its effective domain ${\rm dom}\,H:=\{q\colon\ H(q)<\infty\}$ is the set $\{J(q)<\infty\}$ where 
\begin{equation}
\label{sec.Main.1-00}
 J(q):=
\Pi \big(I_{\{|\ln (1+x)|>1\}}\,e^{-q\ln (1+x)}\big)=\Pi \big(I_{\{|\ln (1+x)|>1\}}(1+x)^{-q}\big). 
\end{equation}
Its interior is the open interval  $(\underline q,\bar q)$ with 
$$
\underline q:=  \inf\{q\le 0\colon\, J(q)<\infty\},\qquad   \bar q:=\sup\{q\ge 0\colon\, J(q)<\infty\}.
$$
Being a   convex function, $H$  is continuous and 
admits finite right and left derivatives on  $(\underline q,\bar q)$.  If $\bar q>0$, then the right derivative satisfies
$$
D^+H(0)=-a_V-\Pi(\bar h(\ln(1+x)))<\infty,
$$ 
though it may be equal to $-\infty $. We do not exclude this case.  

In the formulation of the asymptotic results we shall always assume that $\bar q>0$ and that the equation $H(q)=0$ has a root $\beta\in ]0,\bar q[$. Since $H$
is not constant, such a root is unique.  Clearly, it exists if and only if  $D^+H(0)<0$ and 
$\limsup_{q\uparrow\bar q}H(q)/q>0$. In the case where  $\underline q<0$  the condition  $D^-H(0)>0$
is necessary to ensure that $H(q)<0$ for values $q<0$ which are sufficiently small in absolute value.

 If $J(q)<\infty$, then   
 the process $m=(m_t(q))_{t\ge 0}$ with 
\beq
\label{mart}
m_t(q):=e^{-qV_{t}-tH(q)}
\eeq
is a martingale and
\begin{equation}
\label{sec.Main.1-0}
 \E\, e^{-qV_{t}}= e^{tH(q)}.      
\end{equation}

%\begin{lemma} 
%\label{sup-e}
%Suppose that $Ee^{\e T_1}<\infty$ for some $\e>0$. Let $\beta>0$ be the root of the equation $H(q)=0$. Then 
%\beq
%\label{sup}
%\E\, \sup_{s\le T_1} e^{-\beta V_{s}}<\infty.
%\eeq
%\end{lemma}
%\noindent
%{\sl Proof.}
%Let us take $r>1$ sufficiently close to 1 to ensure the inequality $H(\beta r)<\e$. 
%The process  $e^{-\beta V}=m_s(\beta )$ is a martingale and by the Doob 
%inequality   
%$$
%\E \sup_{s\le T_1} e^{-\beta V_{s}}\le C_r\E   e^{-\beta r V_{s}}=C_r e^{t H(\beta r)}\le C_re^{\e t}.  
%$$
%Therefore, 
%$$
%\E \sup_{s\le T_1} e^{-\beta V_{s}}=\int_0^\infty \E \sup_{s\le t} e^{-qV_{s}} F(dt)\le 
%C_r \E\,e^{\e T_1}<\infty. 
%$$
% The lemma is proven. \fdem 

%\begin{lemma} 
%\label{sup-e}
%Suppose that $Ee^{\e T_1}<\infty$ for some $\e>0$. Let $\beta>0$ be the root of the equation $H(q)=0$. Then for any $q\in ]0,\beta[$
%\beq
%\label{sup}
%\E\, \sup_{s\le T_1} e^{-qV_{s}}<\infty.
%\eeq
%\end{lemma}
%\noindent
%{\sl Proof.}
%Let us take $r\in ]1,\bar q/q[$ sufficiently close to 1 to ensure the inequality $H(qr)-rH(q))<\e$. Then 
%$$
%\E\, m^r_t(q)=e^{t (H(qr)-rH(q))}\le e^{\e t}.  
%$$
%Using Doob's inequality we get that  
%$$
% \E  \sup_{s\le t}m_s(q) \le 1 + \E  \sup_{s\le t}m^r_s(q)\le 1+C_r e^{t (H(qr)-rH(q))}\le 1+C_re^{\e t}.  
%$$
%Since $H(q)<0$,
%$$
%e^{-qV_{s}}  =   m_s(q)e^{s H(q)}\le  \sup_{s\le t}  m_s(q). 
%$$
%Therefore, 
%$$
%\E \sup_{s\le T_1} e^{-qV_{s}}=\int_0^\infty \E \sup_{s\le t} e^{-qV_{s}} F(dt)\le 
%1+C_r \E\,e^{\e T_1}<\infty. 
%$$
% The lemma is proven. \fdem 

\begin{lemma} 
\label{sup-e}
Suppose that $Ee^{\e T_1}<\infty$ for some $\e>0$. Let $\beta\in ]0,\bar q[$ be the root of the equation $H(q)=0$. If $q\in [\beta, \bar q[$ is such that $H(q)\le \e/2$, then   
\beq
\label{sup}
\E\, \sup_{s\le T_1} e^{-qV_{s}}<\infty.
\eeq
\end{lemma}
\noindent
{\sl Proof.}
Let us take $r\in ]1,\bar q/q[$  sufficiently close to 1 to ensure the inequality $H(qr)-rH(q))<(1/2)\e $. Then 
$$
\E\, m^r_t(q)=e^{t (H(qr)-rH(q))}\le e^{(1/2)\e  t}.  
$$
Using Doob's inequality we get that  
$$
 \E  \sup_{s\le t}m_s(q) \le  \E  \sup_{s\le t}m^r_s(q)\le C_r e^{t (H(qr)-rH(q))}\le C_re^{(1/2)\e  t}.  
$$
Note that 
$$
e^{-qV_{s}}  =   m_s(q)e^{s H(q)}\le  e^{tH(q)}\sup_{s\le t}  m_s(q)\le e^{(1/2)\e t}\sup_{s\le t}  m_s(q). 
$$
Therefore, 
$$
\E \sup_{s\le T_1} e^{-qV_{s}}=\int_0^\infty e^{(1/2)\e t}\E  \sup_{s\le t}  m_s(q) F(dt)\le 
C_r \E\,e^{\e T_1}<\infty. 
$$
 The lemma is proven. \fdem 

\begin{corollary} 
\label{coro}
Let $Q:=e^{-V}\cdot P_{T_1}$. If  the conditions of the above lemma are fulfilled and $\E\,|\xi_1|^\beta<\infty$, then $\E\,|Q|^\beta<\infty$. 
%$$
%\E\,\left(\int_{[0,T_1]}e^{-V_t}dt\right)^\beta<\infty. 
%$$
\end{corollary}

{\sl Proof.} Take $q\in ]\beta,\bar q[$ such that $H(q)<\e/2$. Using an obvious estimate and applying the  
H\"older inequality with $r:=q/\beta>1$ and $p:=r/(r-1)$,  
we get that    
$$
\E\,\left(\int_{[0,T_1]}e^{-V_t}dt\right)^\beta\le  \E\, T_1^{\beta} \sup_{s\le T_1} e^{-\beta V_{s}}\le (\E\, T_1^{p\beta})^{1/p}\left(\E\,\sup_{s\le T_1} e^{-q V_{s}}\right)^{1/r} <\infty. 
$$
If $\beta\le 1$, the elementary inequality $||a|+|b||^\beta\le |a|^\beta+|b|^\beta$
implies that 
$$
\E\, |Q|^{\beta}\le c^\beta\E\,\left(\int_{[0,T_1]}e^{-V_t}dt  \right)^\beta+ \E\,e^{-\beta V_{T_1}}\E\,|\xi_1|^\beta <\infty.
$$ 
If $\beta>1$,  we use the inequality $||a|+|b||^\beta\le 2^{\beta-1}(|a|^\beta+|b|^\beta)$ and arrive to the same conclusion. \fdem

\smallskip
Our main result is the following theorem. 

\begin{theorem}
\label{main}
Suppose that there is  $\beta\in ]0,\bar q[$ such that $H(\beta)=0$, 
 $\E|\xi_1|^{\beta}<\infty $. 
 and  $\E\,e^{\e T_1}<\infty$  for some $\e>0$. 
%  $\E\,e^{\e T_1}<\infty$   and $H(\beta)=0$ for some $\beta>0$ such that $\beta\neq \bar q$, the boundary of the effective domain of $H$ (i.e. $H(\beta+\e)<\infty$ for some $\e>0$). Suppose that 
%$\E\,|\xi_1|^{\beta}<\infty $.
If $\sigma\neq 0$ or the law of $|\xi_1|$ has unbounded support,  then 
\beq
\label{maina}
0<\liminf_{u\to \infty}u^{\beta}\Psi(u)\le \limsup_{u\to \infty}u^{\beta}\Psi(u)<\infty. 
\eeq
If $\sigma=0$ and the law of $|\xi_1|$ has bounded support, 
(\ref{maina}) also holds except for the case where  
$0<\Pi(|h|)<\infty$ and  $ \Pi(]-1,0[) \Pi(]0,\infty[)=0$. In the latter case one needs the additional assumption that $\P(T_1\le t)>0$ for any $t>0$. 
\end{theorem}  

{\sl Remark.} In \cite{KPuch} it is shown that in the model where $S$ is a geometric Brownian motion, $\beta:=2a/\sigma-1\le 0$, and $\E\,e^{\e T_1}<\infty$ for some $\e>0$, the ruin probability $\Psi\equiv 1$.

\section{Ruin problem: reduction to a discrete time case}
\label{reduc}
Since we consider the non-life insurance model, ruin may happen only at  jump-times of $P$. This allows us to monitor the process $X^u$ only along the sequence $T_n$, i.e. to reduce the problem to a discrete time setting. For this purpose 
we introduce the discrete time processes $V_{T_n}$ and   
$Y_{n}:=-e^{-V_-}\cdot P_{T_n}$, $n=1,2,...$. Then 
$$
X^u_{T_n}=e^{V_{T_n}}(u- Y_n). 
$$
Note that $Y_n$ is a Markov chain and 
\beq
\label{tY}
Y_{n}=-\sum_{k=1}^n e^{ -V_{T_{k-1}}}\int_{]T_{k-1},T_k]}e^{-(V_{s-}- V_{T_{k-1}})}dP_s
\eeq

Define  the stopping time $\theta^u:=\inf \{n\ge 1\colon Y_n\ge u\}$ with respect to the discrete time 
filtration ${\bf G}=(\cG_n)_{n\ge 1}$ where $\cG_n:=\cF^{R,P}_{T_n}$.  
In the considered model ruin happens only at the jump-times of $P$. Hence,  $\{\tau^u<\infty\}=\{\theta^u<\infty\}$ and  $\Psi (u)= \P(\theta^u<\infty)$.

%The idea of the following lemma is due to Paulsen, see \cite{Paul-93}.  
\begin{lemma}
\label{G-Paulsen} 
If $Y_n\to Y_{\infty}$ almost surely as $n\to \infty$, where  $Y_{\infty}$ is a finite random variable unbounded from above, 
%${\bf H.1}$, ${\bf H.2}$ hold. 
then  for all  $u>0$  
\beq
\label{Paulsen}
\bar G(u)\le\, \Psi(u)
%=\frac{\bar G(u)}{\E\,\bar G(X_{\theta^{u}})\, \vert\, \theta ^{u}<\infty\right]}
\le {\bar G(u)}/{\bar G(0)},
\eeq  
where $\bar G(u):=\P(Y_\infty>u)$. 
%If $\Pi_P((-\infty,0))=0$, then $\Psi(u)=
%\bar G(u)/\bar G(0)$.
\end{lemma}

{\sl Proof.}
Let $\theta$ be an arbitrary stopping time with respect to the  filtration ${\bf G}$. 
As we assume that the finite limit $Y_\infty$ exists,  the random variable 
$$
Y_{\theta,\infty}:=\begin{cases}
-\lim_{N\to \infty } 
\int_{]T_\theta,T_{\theta+N}]}\,e^{-(V_{t-}-V_{T_\theta})}dP_{t},&\theta<\infty, \\
0, & \theta=\infty,
\end{cases} 
$$ 
is well defined.  On the set $\{\theta<\infty\}$
\beq
\label{YX}
Y_{\theta,\infty}=e^{V_{T_\theta}}(Y_{\infty}-Y_\theta)=X_{T_\theta}^u  +e^{V_{T_\theta}}(Y_\infty-u). 
\eeq   Let $\zeta$ be a $\cG_{\theta}$-measurable random variable.  Since the L\'evy process $V$ starts afresh at $\theta$,  the conditional distribution of $Y_{\theta,\infty}$ given $(\theta,\zeta)=(t,x)\in {\mathbb Z}_+\times\bbr$ is the same as the distribution of $Y_{\infty}$. It follows that   
$$
\P(
Y_{\theta,\infty}>\zeta, \ 
\theta<\infty)
=\E\,\bar G(\zeta)\, \Chi_{\{ \theta<\infty\}}.
$$
Thus, if $\P(\theta<\infty)>0$, then
$$
\P(Y_{\theta,\infty}>\zeta, \ 
\theta<\infty)
=\E\,(\bar G(\zeta)\, \vert\, \theta<\infty)\,\P(\theta<\infty).
$$
Noting that  $\Psi(u):=\P(\theta^{u}<\infty)\ge \P(Y_{\infty}>u)>0$,   we deduce from this using (\ref{YX}) that 
\begin{align*}
\bar G(u)&=
\P(
Y_{\infty}>u,\ \tau^{u}<\infty)=
\P(Y_{\theta^u,\infty}>X_{\tau^{u}}^u,\ 
\tau^{u}<\infty)\\
&=\E\,(\bar G(X_{\tau^{u}}^u)\, \vert\, \tau^{u}<\infty)\,\P(\tau^{u}<\infty)\ge \bar G(0)\,\P(\tau^{u}<\infty)
\end{align*}
because $X_{\tau^{u}}^u\le 0$ on   $\{\tau^u<\infty\}$. So, we get the upper  bound 
 in (\ref{Paulsen}). The lower bound is obvious.   
% where  $\Pi_P(]-\infty,0[)=0$, the process $X^u$ crosses zero in a continuous way, i.e.  $X_{\tau^{u}}^u= 0$ on this set. 
 \fdem 
 
In view of the above lemma to prove  Theorem \ref{main} we need to show  that a  finite limit $Y_{\infty}$ exists and is unbounded from above.

\section{The almost sure convergence of $Y_n$ }
\begin{lemma} 
\label{limit}
Suppose that $\E\,e^{\e T_1}<\infty$ and $\E\,|\xi_1|^{\beta\wedge \e\wedge 1}<\infty $ for some $\e>0$ and there is $\beta>0$ such that $H(\beta)=0$. Then the sequence $Y_n$ tends almost surely to a finite random variable $Y_{\infty}$. 
\end{lemma}

{\sl Proof.}
Take  $p\in ]0,\beta\wedge \e\wedge 1[$. According to (\ref{tY})
$$
Y_n-Y_{n-1} =-e^{-V_{T_{n-1}}}\int_{]T_{n-1},T_n]}e^{-(V_{s-}-V_{T_{n-1}})}dP_s=
M_1...M_{n-1}Q_n
$$
where
$$
M_{j}:=e^{-(V_{T_{j}}-V_{T_{j-1}})}, \qquad Q_n=-\int_{]T_{n-1},T_n]}e^{-(V_{s-}-V_{T_{n-1}})}dP_s.
$$
Clearly, $M_1...M_{n-1}Q_n$ is the product of independent random variables, where 
$M_j$ are identically distributed and so are the random variables $Q_n$.  Note that 
$$
\E\, M_1^{p}=\E\, e^{-pV_{T_1}}=\int_0^\infty  \E\, e^{-pV_{t} }    F(dt)=\int_0^\infty  e^{t H(p) } F(dt)< 1. 
$$
Also $\E\, |Q_1|^p<\infty$ (formally we can not use Corollary \ref{coro} because 
we require here a weaker integrability condition on $\xi_1$ but its simple proof need only a minor change). 
. 

%Also for $p\in ]0,1[$
%$$
%\E\, |Q_1|^{p}\le c^p\E\,\left(\int_{[0,T_1]}e^{-V_t}dt  \right)^p+ \E\,e^{-pV_{T_1}}\E\,|\xi_1|^p.
%$$ 
%Applying the H\"older inequality with $p',q'>1$, $1/p'+1/q'=1$, such that $pq'< \beta$ and applying Lemma \ref{sup-e} we get that 
%$$
%\E\,\left(\int_{[0,T_1]}e^{-V_t}dt  \right)^p\le \E\,T_1^p\sup_{t\le T_1}e^{-pV_t} \le 
%\big(\E\,T_1^{pp'}\big)^{1/p'}\big(\E\sup_{t\le T_1}e^{-pq'V_t}\big)^{1/q'}<\infty.
%$$
%So,  
%$\rho:=\E\, M_1^p<1$  and 
%$\E\, |Q_1|^p<\infty$.  
Since   
$\E\, M^p_1...M^p_{j-1}|Q_j|^p=\rho^{j-1}\E\,|Q_1|^p$, we have that 
$$
\E\, \sum_{j\ge 1 } |Y_j-Y_{j-1}|^p<\infty
$$ 
and, therefore, $\sum_{j\ge 1 } |Y_j-Y_{j-1}|^p<\infty$ a.s.. 
Consequently,  $\sum_{j\ge 1} |Y_j-Y_{j-1}|<\infty$ a.s.
and, hence, the sequence $Y_n$ converges a.s. 
to the finite  random variable $Y_\infty:=\sum_{j\ge 1} (Y_j-Y_{j-1})$. \fdem
\smallskip

So, $Y_{\infty}$ is the sum of an absolutely convergent series. Putting $A_0:=1$ and  $A_n:=M_1...M_n$ 
we get that 
$$
Y_\infty=\sum_{n=0}^\infty A_nQ_{n+1}=Q_1+M_1 \sum_{n=1}^\infty  {\frac {A_n}{A_1}}Q_{n+1}  = :  Q_1+M_1 Y_{1,\infty}. 
$$
The random  variable $Y_{1,\infty}$  (we abbreviate here the sum of the series starting from $n=1$)
is independent of $(Q_1,M_1)$ and has the same law as $Y_\infty$.  Using the language of the implicit renewal theory this means that $Y_\infty$ is the solution of the 
distributional (or random) equation  
\begin{equation}
\label{3.1}
Y_{\infty} \stackrel{d}{=}Q+M\,Y_{\infty},\quad
Y_{\infty} \ \ \mbox{independent of}\ (M,Q),
\end{equation}
where $(M,Q):=(M_1,Q_1)$. 

Note that 
$$
Y_\infty=Q_1+A_1Q_2+....+A_{n-1}Q_{n}+ \sum_{k=n}^\infty A_kQ_{k+1}=Q_1+A_1Q_2+....+A_nY_{n,\infty}, 
$$
where 
$$
Y_\infty  \stackrel{d}{=}Y_{n,\infty}:=\sum_{k=n}^\infty \frac {A_{k}}{A_n}Q_{k+1}.  
$$

The following theorem adapted to our needs combines several results of the implicit renewal theory, see Th.  A.6 in \cite{KP2020}.  
\begin{theorem} Suppose that  for  some $\beta>0$,
\begin{align}\label{3.3}
\E\,M^\beta=1, \ \ \ \E\,M^\beta\,(\ln M)^+<\infty, \  \  \ \E\,|Q|^\beta<\infty. 
\end{align}
Then $\limsup  u^\beta \bar G(u)<\infty$.  If  the random variable  $Y_{\infty}$ is unbounded from above, then  
$\liminf  u^\beta \bar G(u)>0$ and  in the case where  the law of $\ln M$ is non-arithmetic, 
$\bar G(u)\sim C_+u^{-\beta}$ where $C_+>0$. 
\end{theorem}
\smallskip

{\sl Proof of Theorem \ref{main}.}  We already know that the sequence of random variables $Y_n$ converges 
a.s. to a finite random variable $Y_\infty$ (see Lemma \ref{limit} whose conditions are weaker than those of Theorem \ref{main}). Moreover, we just proved that $Y_\infty$ is the solution of  distributional equation (\ref{3.1}) and $E\,|Q|^\beta<\infty$ (Corollary \ref{coro}). If  
$Y_\infty$ is unbounded from above  we can use   Lemma \ref{G-Paulsen} relating the asymptotic of the ruin probability with the tail behavior of the distribution function of $Y_\infty$ and get the claimed result from the theorem cited above.  In the next section we show, assuming that  the almost sure limit $Y_{\infty}$ exists,  that $Y_{\infty}$, indeed, is always unbounded,
except the case  where $\sigma=0$, the random variable $\xi_i$ is bounded,  $0<\Pi(|h|)<\infty$ and  $ \Pi(]-1,0[) \Pi(]0,\infty[)=0$ for which we need the  additional assumption  $\P(T_1\le t)>0$ for any $t>0$.  \fdem

%\smallskip 
%\begin{lemma} If the random variables $Q_1$ and $Q_1/M_1$ are unbounded from above, then  
% $Y_{\infty}$ is also unbounded from above. 
%\end{lemma}
%\noindent
%{\sl Proof.}   Since 
%$Q_1/M_1$ is unbounded from above and  independent on $Y_{1,\infty}$, we have  that $\P(Y_{1,\infty}>0)=\P(Y_{\infty}>0)=
%\P(Q_1/M_1+Y_{1,\infty}>0)>0$. 
%  Take arbitrary $u>0$.  Then 
%\bean
%\P(Y_{\infty}> u)&\ge& 
%\P(Q_1+M_1Y_{1,\infty}>u,\, Y_{1,\infty}> 0)
%\ge \P(Q_1>u,\, Y_{1,\infty}> 0)\\
%&=& \P(Q_1>u)\P(Y_{1,\infty}> 0)>0
%\eean
%and the lemma is proven. \fdem
\section{When is $Y_{\infty}$ unbounded from above ?} 
We start from the following elementary consideration assuming that 
the limit $Y_\infty$ exists and is finite. 

\begin{lemma} 
\label{tail}
If the random variables $Q_1$ and $Y_n/A_n$ for some $n\ge 1$ 
%$Q_1/A_2+Q_2/A_1$ 
are unbounded from above,  then  $Y_{\infty}$ is  unbounded from above. 
\end{lemma}
\noindent
{\sl Proof.}  Take $N$ such that $\P(Y_{\infty}>-N)>0$. Since 
$Y_n/A_n$ is unbounded from above and independent of the random variable $Y_{n,\infty}$ which has the same law as $Y_\infty$, we get that 
\bean
\P(Y_\infty>0)&=&
\P(Y_{n}+ A_{n} Y_{n,\infty}>0)= \P(Y_n/A_n+Y_{n,\infty}>0)\\
&\ge& \P(Y_n/A_n\ge N, \,Y_{n,\infty}>-N)\\
&=&\P(Y_n/A_n\ge N)\P(Y_{n,\infty}>-N)>0.
\eean 
Thus, $\P(Y_{n,\infty}>0)=\P(Y_{\infty}>0)>0$. 
%So,  $\P(Y_{k,\infty}>0)=\P(Y_{\infty}>0)>0$ for all $k\ge 1$. 
Furthermore, if the random variable  $Q_1$ is unbounded from above and $\P(Y_{\infty}>0)>0$,  then 
$Y_\infty$ is also unbounded from above. Indeed, 
take arbitrary $u>0$, then 
 \bean
\P(Y_{\infty}> u)&\ge& 
\P(Q_1+A_1 Y_{1,\infty}>u,\, Y_{1,\infty}> 0)
\ge \P(Q_1>u ,\, Y_{1,\infty}> 0)\\
&=& \P(Q_1>u )\P(Y_{1,\infty}> 0) = \P(Q_1>u )\P(Y_{\infty}> 0)>0. 
\eean
%\bean
%\P(Y_{\infty}> u)&\ge& 
%\P(Y_n+A_n Y_{n,\infty}>u,\, Y_{n,\infty}> 0)
%\ge \P(Y_n>u ,\, Y_{n,\infty}> 0)\\
%&=& \P(Y_n>u )\P(Y_{n,\infty}> 0)>0
%\eean
The lemma is proven. \fdem

\smallskip
We shall use the above lemma with $n=1$ when $Y_1/A_1=Q_1/A_1=Q_1/M_1$ or with 
$n=2$ when $Y_2/A_2=Q_1/A_2+Q_2/M_2$. 

\smallskip
The  arguments below  use the following  observation.  Let $\zeta$ be a real-valued random variable and let  $\eta$ be a random variable with values in a Polish space.  Let  $\P_\zeta$ and $\P_\eta$ be their distributions and 
let $\P_{\zeta|x}$ be a regular conditional distribution of $\zeta$ given $\eta=x$. If for any real $N$ the set  $\{x \colon \P_{\zeta | x} ( [N,\infty))>0\}$ is not a $\P_\eta$-null set, then $\zeta$ is unbounded from above.

\begin{lemma} 
\label{Wiener}
Let $K>0$,   $\sigma\neq 0$ and $t>s>0$. Then the  random variables 
\beq
\zeta:= Ke^{\sigma W_t} -\int_0^t e^{\sigma W_r}dr, \qquad  \tilde\zeta:=Ke^{\sigma (W_t-W_s)} -e^{\sigma W_t}\int_0^t e^{\sigma W_r}dr, 
\eeq
are unbounded from above.
% (the sign of $\sigma$ does not matter  since $-W$ is again a Wiener process)). 
\end{lemma}
{\sl Proof.}  Recall that the Wiener measure charges any open ball   in the space $C_0([0,T])$ of continuous functions $x_.=(x_t)_{t\in [0,T]}$ with $x_0=0$ (in other words, it is of full support on this space).  Let $g$ be a continuous function on 
$C_0([0,T])$ such that $g(x^0_.)=N$ and let $\e>0$.  Then there is a $\delta>0$ such that  
$|g(x_.)-g(x^0_.)|<\e$ on the open ball $\{x_.\colon |x_.-x^0_.|<\delta\}$. In particular, $g(x_.)> N-\e$ on this ball. 
With this, the claims of the lemma are obvious. 
\fdem 

%
%Without loss of generality we may assume that $t=1$.  Recall that the conditional 
%law of $W=(W_s)_{s\le 1}$ given $\eta:=W_1=x$ is the law of the Brownian 
%bridge $B^x=(B^x_s)_{s\in 1}$  with $B^x_0=0$ and $B^x_1=x$, which coincides with the
%law of process $(W_s+(x-W_1)s)_{s\le  t}$. The law of the random variable $\int_0^1 e^{W_s+(x-W_1)s}ds$ charges any interval  $]0,\e]$.    Therefore, 
%$
%\P(\P_{\zeta | x} (\zeta\ge N)>0)>0$  for sufficiently large $x$ (e.g., when  $ K e^{\sigma x}-1\ge N$), i.e.  
%unbounded from above. 
%
%The conditional law of $\tilde \zeta$ given  $\eta:=W_1=x$ is the unconditional law of the random variable 
%$$
%e^{\sigma x} \Big (K e^{-\sigma B^x_s} - \int_0^1 e^{B^x_r}dr\Big) \stackrel{d}{=}
%e^{\sigma x} \Big (K e^{-\sigma (W_s+(x-W_1)s)} - \int_0^1 e^{\sigma(W_r+(x-W_1)r)}dr\Big)
%$$
% \fdem

\smallskip
Note that in the above lemma the sign of $\sigma$ does not matter  since $-W$ is again a Wiener process. 

%Using this observation with get that the random variable  
%\beq
%\zeta:= Ke^{\kappa W_t} -\int_0^t e^{\kappa W_s}ds, 
%\eeq
%where  $K>0$ and $\kappa \neq 0$ (the sign of $\kappa$ does not matter  since $-W$ is again a Wiener process),  is unbounded from above.
%% Indeed, changing the variable in the integral and the  parameter $\kappa$ we may assume  without loss of generality that $t=1$. 
% 
% The conditional 
%law of $W=(W_s)_{s\le t}$ given $\eta:=W_t=x$ is the law of the Brownian 
%bridge $B^x=(B^x_s)_{s\in t}$  with $B^x_0=0$ and $B^x_t=x$, which coincides with the
%law of process $(W_s+(x-W_t)s/t)_{s\le  t}$. The law of the integral $\int_0^t e^{W_s+(x-W_t)s/t}ds$ charges any interval  $]0,\e]$.    Thus, 
%$
%\P(\P_{\zeta | x} (\zeta\ge N)>0)>0$  for sufficiently large $x$ (e.g., when  $ e^{\kappa x}-1\ge N$). 
%
%In a similar way we can obtain the same property of the random variable
%\beq
%\label{z2}
%\zeta:=K_1e^{\kappa (W_t-W_s)} -K_2e^{\kappa W_t}\int_0^t e^{\kappa W_r}dr. 
%\eeq
%Indeed,  given 

Let $\bar V:=V-\sigma W$. Then  
$$
Q_1=e^{-V_{T_1}}|\xi_1| -c\int_0^{T_1} e^{-V_r}dr\ge |\xi_1|e^{-\bar V_{T_1}}e^{-\sigma W_{T_1}}-c
\sup_{r\le T_1} e^{-\bar V_r} \int_0^{T_1}e^{-\sigma W_r}dr. 
$$
Taking into account the independence  of $\xi_1$, $T_1$, $\bar V$ and $W$, we obtain using conditioning and Lemma \ref{Wiener}, that $Q_1$ is unbounded from above. 

%Using conditioning with respect to $\xi_1$, $T_1$ and $V-\sigma W$ we get that for the case $\sigma\neq 0$ the random variable
%$$
%Q_1=e^{-V_{T_1}}|\xi_1| -c\int_0^{T_1} e^{-V_s}ds
%$$
% is unbounded from above.  
 
 The random variable $Q_1/M_1$ is unbounded from above  if and  only if $|\xi_1|$ has this property. 
 To cover the general case we consider the random variable  
%Inspecting the expressions  
%$$
%Q_1=-c\int_0^{T_1} e^{-V_s}ds+e^{-V_{T_1}}|\xi_1|, \qquad  Q_1/M_1=-c\int_0^{T_1} e^{V_{T_1}-V_s}ds+|\xi_1|.
%$$ 
%we see that in the case where $\sigma\neq 0$, then  $Q_1$ is unbounded from above and so is 
%$Q_1/M_1$, if $|\xi_1|$ is unbounded. 
\begin{align}
Y_2/A_2
%=Q_1/A_2+Q_2/A_1
&=e^{V_{T_2}-V_{T_1}}|\xi_1|+ |\xi_2|-ce^{V_{T_2}}\int_{0}^{T_2}e^{-V_r}dr \notag\\
&\ge  e^{V_{T_2}-V_{T_1}}|\xi_1|-ce^{V_{T_2}}\int_{0}^{T_2}e^{-V_r}dr.   \label{eq:720}
\end{align}
The same arguments as above allow us to conclude that the random variable on the right-hand side is 
unbounded from above.

\smallskip
Thus, if  $\sigma\neq 0$,  the random variable  $Y_\infty$ is unbounded without further assumptions. 

\smallskip
It  remains to consider the case
$\sigma=0$.  The further arguments  are based on decompositions of  L\'evy processes without Gaussian component  into sums of  independent L\'evy processes and the fact that a Poisson process has on a finite interval any  number of jumps with strictly positive probability.    The  reasoning depends on the behavior of the L\'evy measure.

{\bf 1.} Suppose that $\Pi(]-1,0[)>0$ and, hence,  $\Pi(]-1, -\e[)>0$ for some $\e>0$.   
We decompose $V$ into a sum of two independent processes  $V=V^{(1)}+V^{(2)}$ where $V^{(2)}=V-V^{(1)}$ and  
$$
V^{(1)}_s:=
I_{]-1, -\e[}\ln (1+x)*\mu_s +I_{]\e,\infty [}\ln (1+x)*\mu_s. 
$$ 
Note that $V^{(1)}$ is the sum of two independent compound Poisson processes with negative and positive jumps respectively. 
%V^{(2)}_s&=& as- I_{]-1, -\e[}h*\nu_s+I_{]-\e, 0[}h*(\mu-\nu)_s+I_{]-\e,0[}(\ln (1+x)-h)*\mu_s. 
Using conditioning with respect to the random variables $V^{(2)}$, $T_1$, $\xi_1$, which are independent of $V^{(1)}$, we easily obtain that 
$Q_1$ is unbounded from above because for any $K>0$ the random variable 
\beq
\zeta:= Ke^{-V^{(1)}_t} -\int_0^t e^{-V^{(1)}_r}dr, \quad K>0, 
\eeq
is unbounded from above. The latter property is clear: we may consider trajectories where $V^{(1)}$ has only  negative jumps and these jumps  are concentrated
near $t$ leading to large values for the first term and a small impact on the value of the integral. 

If $|\xi_1|$ is unbounded, then $Q_1/M_1$ is unbounded from above and, hence, so is $Y_\infty$.

If $|\xi_1|$ is bounded and $\Pi(]0,\infty[)>0$ (hence for some sufficiently small $\e>0$ both components of $V^{(1)}$ are nontrivial compound Poisson processes),   we can check that 
  $Y_2/A_2$ is unbounded from above.  Indeed, from \eqref{eq:720},
 $$
 Y_2/A_2\ge  
% e^{V_{T_2}}\Big (e^{-V_{T_1}}|\xi_1|-c\int_{0}^{T_2}e^{-V_r}dr\Big)= 
 e^{V_{T_2}-V_{T_1}}\Big (|\xi_1|-ce^{V_{T_1}}\int_{0}^{T_1}e^{-V_r}dr-c \int_{T_1}^{T_2}e^{-(V_r-V_{T_1})}dr\Big).
 $$
 Using conditioning, we reduce the problem to checking that  for any $t>s>0$ and $K>0$ the random  
 variable
 \beq
 \label{V1} 
%e^{V^{(1)}_t }\Big (K e^{-V^{(1)}_s}-\int_{0}^{t}e^{-V^{(1)}_r}dr\Big){\color{red}=
e^{V^{(1)}_t -V^{(1)}_s }\Big (K -e^{V^{(1)}_s}\int_{0}^{s} e^{-V^{(1)}_r}dr-\int_{s}^{t}e^{-(V^{(1)}_r-V^{(1)}_s)}dr\Big)
\eeq
is unbounded from above.  The compound Poisson process $(V^{(1)}_r)_{r\in [0,s]}$ has  positive and negative jumps with intensities bounded away from zero and, therefore,  for any  $\e'>0$  there is a non-null set $\Gamma_1$ on which  
$$
e^{V^{(1)}_s}\le 1, \qquad \int_{0}^{s} e^{-V^{(1)}_r}dr\le \e'. 
$$   
The  process $(V^{(1)}_r-V^{(1)}_s)_{r\in [s,t]}$ is  independent of $(V^{(1)}_r)_{r\in [0,s]}$ and has the same jump properties. For an $\e'>0$ and $N>0$ there is a non-null set $\Gamma_2$ on which 
$$
e^{V^{(1)}_t -V^{(1)}_s } \ge N, \qquad \int_{s}^{t}e^{-(V^{(1)}_r-V^{(1)}_s)}dr \le \e'.
$$
Taking $\e'\le K/4$ we get that on the non-null set $\Gamma_1\cap \Gamma_2$ the random variable defined in (\ref{V1}) is larger than  $NK/2$ and, therefore, is unbounded from above. 

%As a consequence
%$$ e^{V^{(1)}_s}\int_{0}^{s} e^{-V^{(1)}_r}dr \le \eN. $$
%The paths of the process $(V^{(1)}_r-V^{(1)}_s)_{r\in [s,t]}$ independent of $(V^{(1)}_r)_{r\in [0,s]}$ have the same jump properties, and therefore on a non-null set $\Gamma_2$ additionally the following inequalities hold:
% With this, it is easily seen that the random variable defined in (\ref{V1}) 
%is unbounded from above. 

%
%But  this is rather obvious. Indeed, we note that $V^{(1)}$ has, with a strictly positive probability,
% only negative jumps on $[0,s]$, only positive ones on $]s,t]$, and a number of sufficiently big jumps near ${\color{red} {t}}$. Therefore, the large value of  $e^{V^{(1)}_t }$ jointly with the strict positivity of the expression in the parentheses ensure this to be the case. 

\smallskip
{\bf 2.}  Suppose that $\Pi(]0,\infty[)=0$ and $\Pi(|h|)=\infty$.
Then for any $N>0$ 
%and $t>0$ (how does the t come in here?)\\
there exists $\e>0$ such that $ \tilde R_\e:=\Pi(I_{]-1,-\e[}|h|)\ge N$. Let us consider the decomposition   $V=V^{(1)}+V^{(2)}$ into a sum of two independent processes, this time with   
\bean
V^{(1)}_s: = I_{]-1,-\e]}h*(\mu -\nu)_s+
I_{]-1,-\e]}(\ln (1+x)-h)*\mu_s.
%V^{(4)}_s&=& as- I_{]0, \e]}h*(\mu-\nu)_s+ I_{]0, \e]}(\ln (1+x)-h)*\mu_s. 
\eean
Again, even if $|\xi_1|$ is bounded, the random variable  $Y_2/A_2$ is unbounded from above. 
Indeed, on the non-null set where the process $(V^{(1)}_r-V^{(1)}_s)_{r\in [s,t]}$  has no jumps and 
$$
e^{V^{(1)}_s}\int_{0}^{s} e^{-V^{(1)}_r}dr\le K/2
$$  
the random variable given by (\ref{V1})  dominates  the value
$$
e^{\tilde R_\e(t-s)}\Big( K/2 -\int_s^t e^{\tilde R_\e (s -r) } dr \Big)\ge e^{N(t-s)}( K/2 -1/N).   
$$
This implies the desired property.

\smallskip
The remaining subcase in which $\Pi(]-1,0[)>0$ but $\Pi(]0,\infty[)=0$ will be considered later. 

\smallskip
{\bf 3.}  Suppose that $\Pi(]-1,0[)=0$ and $ \Pi(h)=\infty$.  Then for any $N>0$ 
%and $t>0$\\ (again how does t come in here?)\\
there exists $\e>0$ such that $ R_\e:=\Pi(I_{]\e,1]}h)\ge N$. 
We consider again a decomposition  $V=V^{(1)}+V^{(2)}$ in a sum of two independent processes. The process  $V^{(1)}$ is now defined as
\bean
V^{(1)}_s&:=&I_{]\e, 1]}h*(\mu -\nu)_s+
% I_{]-1, -\e[}h*\mu_t+I_{]-1, -\e[}(\ln (1+x)-h)*\mu_t
I_{]\e, 1]}(\ln (1+x)-h)*\mu_s.
%V^{(4)}_s&=& as- I_{]0, \e]}h*(\mu-\nu)_s+ I_{]0, \e]}(\ln (1+x)-h)*\mu_s.
\eean
The set on which $V^{(1)}$ has no jumps on $[0,t]$ has strictly positive probability. On this set we derive the bound
\beq
Ke^{-V^{(1)}_t} -\int_0^t e^{-V^{(1)}_r}dr\ge Ke^{R_\e t}-\frac 1{R_\e} e^{R_\e t} \ge e^{Nt}(K-1/N)
\eeq
and we conclude as above that $Q_1$  is unbounded from above. 

If $|\xi_1|$ is unbounded, then $Q_1/M_1$ is unbounded from above as well as $Y_\infty$. 

If $|\xi_1|$ is bounded we again can prove that  $Y_2/A_2$ is unbounded from above. Indeed, on the set on which the process  $V^{(1)}$ has no jumps on $[0,s]$ (hence, it is growing linearly on this interval) we conclude that 
\beq
e^{V^{(1)}_t }\Big (K e^{-V^{(1)}_s}-\int_{0}^{t}e^{-V^{(1)}_r}dr\Big)\ge  e^{V^{(1)}_t -V^{(1)}_s}\Big( K-\frac 1{R_\e}  -\int_s^t e^{V^{(1)}_s -V^{(1)}_r } dr \Big).  
 \eeq
Let us define $J=I_{]\e, 1]}\ln(1+x)*\mu$, then  $V^{(1)}_t -V^{(1)}_s=-R_\e(t-s)+J_t-J_s$ and  $V^{(1)}_s -V^{(1)}_r\le R_\e(t-s)- (J_r-J_s)$. Since the increment $(J_r-J_s)$ may have arbitrary many jumps (all of a size larger than $\ln(1+\e)$) and this happens independently of $J_{r\in [0,s]}$, we get that 
the  random variable given by (\ref{V1}) is unbounded from above.

\begin{lemma} Suppose that $\sigma=0$, $\Pi(]-1,0[)=0$ and $ 0<\Pi(h)<\infty$. If $\P(T_1\le t)>0$ for any $t>0$, then $Y_\infty$ is unbounded from above. 
\end{lemma}
{\sl Proof.} In this case $V_t=-a_ht+L_t$ where  the process   $L_t:=\ln (1+x)*\mu_t$ is increasing. The constant $a_h:=\Pi(h)-a$ is strictly positive because  the equality 
$\ln \E e^{-\beta V_1}=0$  with $\beta>0$ may hold only if 
%Obviously, 
%$$Y_{\infty}\ge a^_P\int_0^Te^{-V_t}dt =a^0_P\int_0^Te^{-V_t}dt.$$  
$\P(V_1<0)>0$.   
We have 
\beq
\label{Yn}
Y_n:=-e^{-V_-}\cdot P_{T_n}=\sum_{k=1}^n e^{a_h T_k-L_{T_k}}|\xi_k|-c\int_{0}^{T_n} e^{a_hs-L_s}ds. 
\eeq

Let $\delta>0$. Put $L^\delta_t:=I_{\{x<\delta\}}\ln (1+x)*\mu_t$. Then for every $t\ge 0$
$$
\E\, L^\delta_t=
t \Pi(I_{\{x<\delta\}}\ln (1+x))\downarrow 0,  \quad \delta \downarrow 0.
$$ 
Thus, for a sufficiently small $\delta$ the set $\{L^\delta_t\le 1\}$ is non-null. Also 
non-null is the set $\{L_t-L^\delta_t=0\}$.  Their intersection, due to the independence of $L^\delta_t$ and $L_t-L^\delta_t$, is non-null and  the larger set $\{L_t\le 1\}$ is also non-null. 

Take  $K', \kappa>0$ such that $\P(1/K'\le T_1\le K')>0$ and 
$\P(|\xi_1|>\kappa)>0$.  Then    $\P(\inf_{k\le n}|\xi_k|>\kappa)>0$ for each $n$ due to the independence of the random variables $\xi_k$. 
%Put $p:=e^{a_h/K}$. 
%Fix  $p>1$ and $K>0$ such that $\P(e^{a_hT_1}>p,\; T_1\le K)>0$. 

The independence of the interarrival times and the assumption on their distribution imply that for any $n$ the set 
$$
\{1/K'\le T_1\le K',\; T_k-T_{k-1}\le 2^{-(k-1)}, \; 1< k \le n\}
$$ 
is non-null  as well as its intersection with the set $\{\inf_{k\le n}|\xi_k|>\kappa, L_{K'+2}\le 1\}$. On this intersection, which we denote by $\Gamma_n$,  we minorate the sum in the right-hand side of (\ref{Yn}) by replacing  $L_{T_k}$ by unit and majorate the integral by   replacing  $L_{s}$  by zero.  Taking into account that $T_n\le T_1+1\le K'+1$, we get in this way  that  
$$
Y_n\ge  \frac{\kappa}e \sum_{k=1}^n e^{a_h T_k} -\frac {c}{a_h}e^{a_h(T_1+1)} \ge \Big(\frac{\kappa}e n - \frac {c}{a_h}e^{a_h}\Big)e^{a_hT_1} \ge \Big(\frac{\kappa}e n - \frac {c}{a_h}e^{a_h}\Big)e^{a_h/K'}
$$
where the last inequality holds for $n\ge c e^{a_h+1}/(\kappa a_h)$. 
Thus, for any fixed $N$ we have on the set $\Gamma_n$  the bound  $Y_n\ge N$ when $n$ is sufficiently large. 

Take $y>0$ such that $\P(Y_{\infty}>-y)>0$. Due to the independence the set   
$$
\Gamma_n\cap \{T_{n+1}-T_n\le 1,\, Y_{n+1,\infty}>-y\}
$$ 
is  non-null  and on this set for sufficiently large $n$ we have that 
\bean
Y_{\infty}&=&Y_{n+1}+A_{n+1} Y_{n+1,\infty}=Y_n+A_n Q_{n+1}+A_{n+1} Y_{n+1,\infty}\\
&=& Y_n
+ e^{a_hT_{n+1}-L_{T_{n+1}}}|\xi_{n+1}| -c\int_{T_{n}}^{T_{n+1}}e^{a_hs-L_s}ds + 
e^{a_hT_{n+1}-L_{T_{n+1}}}Y_{n+1,\infty}\\
&\ge & Y_n -(c/a_h)e^{a_h(K'+2)} - e^{a_h (K'+2)}y\ge N  -(c/a_h+y)e^{a_h(K'+2)}. 
\eean
Since $N$ is arbitrary, the random variable $Y_{\infty}$ is unbounded from above. \fdem 

\smallskip
\begin{lemma} Suppose that $\sigma=0$, $ \Pi(]0,\infty[)=0$, and $0<\Pi(|h|)<\infty$. If $\P(T_1\le t)>0$ for any $t>0$, then $Y_\infty$ is unbounded from above. 
\end{lemma}
{\sl Proof.}   In this case $V_t=b_ht-L_t$ where the process  $L_t:=-\ln (1+x)*\mu_t$ is increasing and the  $b_h:=a-\Pi(h)>0$.  
We have 
$$
Y_n:=-e^{-V_-}\cdot P_{T_n}=-\sum_{k=1}^n e^{-b_h T_k+L_{T_k}}\xi_k-c\int_{0}^{T_n} e^{-b_hs+L_s}ds. 
$$
As in the proof of the previous lemma, we get that 
for any $t\ge 0$ the set $\{L_t\le 1\}$ is a non-null set,  there are  constants $\kappa$ and $K'$ such that 
 $\P(|\xi_k|>\kappa)>0$ for each $k$ and  
 $\P(1/K'\le T_1\le K')>0$. 
For any $n$ the set 
$$
\{1/K'\le T_1\le K',\; T_k-T_{k-1}\le 2^{-(k-1)}, \; 1< k \le n\}
$$ 
is non-null  and its intersection with the set $\{\inf_{k\le n}|\xi_k|>\kappa, L_{K'+2}\le 1\}$ is non-null. On this intersection, denoted $\Gamma_{n}$,  
$$
Y_n\ge \kappa e^{-b_hT_{1}}\sum_{k=1}^n e^{-b_h (T_k-T_1)} -\frac {ce}{b_h} \ge \kappa e^{-b_h(K'+1)}n - \frac {ce}{b_h}\ge  N
$$
for any $N$ when $n$ is sufficiently large. 

Take $y>0$ such that $\P(Y_{\infty}>-y)>0$. Then the set 
$$
\Gamma_{n} \cap\{T_{n+1}-T_n\le 1,\;Y_{n+1,\infty}>- y\}
$$
is non-null and on this set for sufficiently large $n$ we have that 
\bean
Y_{\infty}&=&Y_n+ e^{-b_hT_{n+1}+L_{T_{n+1}}}|\xi_{n+1}| -c\int_{T_{n}}^{T_{n+1}}e^{-b_hs+L_s}ds + 
e^{-b_hT_{n+1}+L_{T_{n+1}}}Y_{n+1,\infty}\\
&>& Y_n -ce - ey\ge N -e(c+y). 
\eean
Since $N$ is arbitrary, the random variable $Y_{\infty}$ is unbounded 
from above. \fdem 

%\bigskip
%\noindent
%{\bf References}

 \end{document}